\documentclass[a4paper,12pt]{article}
\usepackage{amsfonts,amssymb,xcolor}
\title{ Locally convex topological algebras  of  generalized functions: compactness and nuclearity in a nonlinear context}

\author{J. Aragona,  S. O. Juriaans,\\
Instituto de Matem\'atica e Estat\'istica, Universidade  de S\~ao Paulo.\\
 J. F. Colombeau$^{+}$,\\
 Institut Fourier, Universit\'e de Grenoble.}

\date{}

\begin{document}
\maketitle
\textit{Dedicated to the memory of Leopoldo Nachbin.}\\
\\

\begin {abstract} In this paper we introduce Hausdorff locally convex algebra topologies on subalgebras of the whole algebra of nonlinear generalized functions. These topologies are strong duals of Fr\'echet-Schwartz space topologies and even strong duals of nuclear Fr\'echet space  topologies. In particular any bounded set is relatively compact and one benefits from all deep properties of nuclearity. These algebras of generalized functions contain most of the classical irregular functions and distributions. They are obtained by replacing the mathematical tool of $\mathcal{C}^\infty$ functions in the original version of nonlinear generalized functions by the far more evolved tool of holomorphic functions. This paper continues the nonlinear theory of generalized functions in which such locally convex topological  properties were strongly lacking up to now.\\
\end{abstract}


\textit{(+) the work of this author was done under financial support of FAPESP, processo 2011/12532-1, and thanks to the hospitality of the IME-USP.}\\
\textit{(+) corresponding author jfcolombeau@ime.usp.br}\\

\textbf{1. Introduction}. 
We denote by $\mathcal{G}(\Omega)$ the special (or simplified) algebra of nonlinear generalized functions on $\Omega$. We construct subalgebras of $\mathcal{G}(\Omega)$ that enjoy Hausdorff locally convex topologies  suitable for the development of a  functional analysis: they have very good topological properties in particular concerning compactness (strong duals of Fr\'echet-Schwartz spaces) at the same time as they are compatible both with partial derivatives and nonlinearity. These properties even extend to nuclearity. Nuclear spaces \cite{Pietsch,Schaefer,Yosida} are an extension of finite dimensional spaces in which finite sums are replaced by convergent series in a natural topology: a difference with Hilbert spaces lies in that in the nuclear (DFN) algebras  the bounded sets are relatively compact (even much more: they resemble to some extent the finite dimensional bounded sets) at the price of the use of nonmetrizable topologies which are topological inductive limits of a sequence of separable Hilbert spaces linked by nuclear inclusion maps.\\

 At the beginning of the theory of nonlinear generalized functions, initiated inside the group of Prof. Leopoldo Nachbin \cite{Colombeau1, Colombeau2, AraBia, Bialivre}, a Hausdorff topology  was  introduced  by H.A. Biagioni in \cite{Bialivre} p. 44-45.  Later it was rediscovered  independently by D. Scarpalezos in  \cite{Scarpatopo}. It was  called "sharp topology" by D. Scarpalezos 
\cite{Nedel, Scarpatopo, Scarpatopo2} who gave the impulse for its use. The sharp (or Scarpalezos) topology was at the origin of  numerous works \cite{Bialivre, Scarpatopo, Scarpatopo2, AragonaJu,   AragonaFerJu1,  AragonaFerJu2,  AragonaFerJu3, Delcroix,  Garetto, Garetto2, GarettoV, AragonaJuriaans, AragonaFernandez,  AragonaGarcia, Oberv}.  Besides its good properties  it is not a usual Hausdorff locally convex algebra topology on the field of real or complex numbers as this would have been welcome for a development of functional analysis in the context of nonlinear generalized functions in order to benefit of compactness and of the well developed mathematical  theory of locally convex spaces.   We should indicate however that in  \cite{Novo, AragonaFerJu2, Oberv} a theory is started leading to a notion of compactness in the sharp topology \cite[Theo. 2.12] {Oberv}.\\

Therefore one could search for a Hausdorff locally convex topology with  compactness properties reminiscent  of the topologies of distribution theory. We do not obtain such a topology on the whole algebra of nonlinear generalized functions, but on  subalgebras of $\mathcal{G}(\Omega)$ that anyway contain most distributions. Such subalgebras are constructed in section 2 while the  topology is constructed in section 4, from a basic compactness result proved in section 3. The far more evolved nuclear topologies are constructed in section 7 from  properties of the theory of nuclear maps and nuclear spaces. These topological subalgebras are equipped with  Hausdorff locally convex topologies of strong duals of Fr\'echet-Schwartz spaces and strong duals of nuclear Fr\'echet spaces.  In the same way as the original version of nonlinear generalized functions stemmed from the canonical Hamiltonian formalism \cite{Colombeau1,Colombeau4, ColombeauGsponer} which suggested to extend the calculations on $\mathcal{C}^\infty$ functions to irregular functions and distributions, the intuitive physical ideas at the origin of the present work stemmed from dimensional renormalization and the  observations of nonlocality of space-time which suggested the use of analytic functions instead of $\mathcal{C}^\infty$ functions. The intuitive mathematical ideas stemmed from the remark that in the context of nonlinear generalized functions analytic representatives appear sufficient to represent most of the irregular functions and distributions of physics, therefore one could restrain the theory to such representatives so as to benefit from the remarkable mathematical properties of analytic functions, even at the price of far more difficult proofs. Indeed the present context can be considered as the restriction of the context used so far to analytic representatives; then it benefits of very different and better properties. \\

%\end{document}

\textbf{2. Construction of a  subalgebra $\mathcal{G}_a(\Omega)$}. 
 We set $z=x+iy, x\in\Omega\subset\mathbb{R}^k$ which is assumed to be open and connected, $y\in\mathbb{R}^k, \zeta=\xi+i\eta, \ \xi,\eta \in \mathbb{R}$.\\

We use as auxiliary ingredient  a family of nonvoid open sets $O_n \subset \Omega$ such that $\forall n  \ O_{n+1}\subset O_n$ and $d(O_{n+1},\mathcal{C}O_n)>0$
if $d(O_{n+1},\mathcal{C}O_n)$ is the distance between $O_{n+1}$ and the complement of $O_n$ in $\Omega$.
Here are various examples:\\
\textit{i)} $ x_0\in \Omega$ and $O_n=\{x\in \Omega / d(x,x_0)<\frac{1}{n}\}$,\\
\textit{ii)} $ x_0\in \partial\Omega$ and $O_n=\{x\in \Omega / d(x,x_0)<\frac{1}{n}\}$,\\
\textit{iii)} $O_n=\{x\in \Omega / |x|>n\}$ when $\Omega$ unbounded,\\
\textit{iv)} $O_n=\{x\in \Omega / d(x,\mathcal{C}\Omega)<\frac{1}{n}\}$.\\

A weaker concept could be used since we basically use only an accumulation point in each $O_n$. The choice of such a family of sets will appear clear in the proof of theorem 1 where they play the basic role. 
The construction below of the topological subalgebras of  $\mathcal{G}(\Omega)$ depends on such a sequence  $(O_n)_n $, supposed to be given from now on.\\

  We set
\begin{equation} A_n:=\{(z,\zeta)\in \mathbb{C}^k\times\mathbb{C}  \   / \ x\in\Omega, x\not\in O_{n+1}, \ |arg \zeta|<\frac{1}{n}, \ 0<|\zeta|< \frac{1}{n}, |y_i|<\frac{1}{n} |\zeta|, 1\leq i\leq k\},\end{equation}

\begin{equation} B_n:=\{(z,\zeta)\in \mathbb{C}^k\times\mathbb{C} \   / \  x\in O_n, \ 0<|\zeta|< \frac{1}{n}, |y_i|<\frac{1}{n},\ 1\leq i\leq k\},\end{equation}

\begin{equation} V_n=  A_n\cup B_n.\end{equation}
\\
The definition of $V_n$ has been chosen since it is convenient for the sequel. $V_n$ is a connected open set in $\Omega\times \mathbb{R}^k\times\mathbb{C}$ (variables $(x,y,\zeta))$, having $\Omega\times \{0_{\mathbb{R}^k}\}\times\{0_{\mathbb{C}}\}$ (variables $(x,y=0,\zeta=0)$ on its boundary. The family $\{V_n\}$ is decreasing and  $\bigcap V_n=\oslash$. The choice of these sets will appear clear from  examples such as the inclusion of the distributions with compact support in the algebra of generalized functions.  Let $\phi:\Omega \longmapsto\mathbb{R}^*_+$ be a strictly positive continuous function. We consider the vector space

\begin{equation} H_{n,\phi}:=\{f:V_n\longmapsto \mathbb{C} \ holomorphic \ / \ \exists \ \ const>0 \  / |f(z,\zeta)|\leq const.\frac{\phi(x)}{|\zeta|^n} \ \forall (z,\zeta) \in V_n\} .\end{equation}

We equip $H_{n,\phi}$ with the norm 

\begin{equation} \|f\|_{n,\phi} :=sup_{(z,\zeta)\in V_n} |\zeta|^n \frac{ |f(z,\zeta)|}{\phi(x)}.\end{equation}
\\
Therefore $|f(z,\zeta)|\leq\|f\|_{n,\phi}\frac{\phi(x)}{|\zeta^n|}$  and  $\|f\|_{n,\phi}$ is the smallest value of $const$ in (4). The normed space
$(H_{n,\phi},\| \ \|_{n,\phi}) $  is a Banach space from the classical properties of holomorphic functions. If $n'\geq n$ and $\psi\geq \phi$, then $H_{n,\phi} \subset H_{n',\psi}$ with the inclusion map $f\longmapsto f_{|V_{n'}}$ which is injective by analyticity since $\Omega$, therefore $V_n$,  is connected. Further it is obvious that
 $$H_{n,\phi}.H_{n,\phi}\subset H_{2n,\phi^2} $$
 where  .  stands for the classical product of functions.\\

One defines an algebra $\mathcal{G}_a(\Omega)$ as the algebraic  inductive limit of the vector spaces $H_{n,\phi}$, when $n\rightarrow +\infty$ and when $\phi$ ranges over all continuous strictly positive  functions $\Omega\longmapsto\mathbb{R}^*_+$, ordered by the natural order from $\mathbb{R}^*_+$. From Cauchy's inequality for the derivation in $z$-variable it is clear that $\mathcal{G}_a(\Omega)$ is a differential algebra; indeed the last inequality in (1) forces to apply Cauchy's formula on a disk of radius $\frac{const}{|\zeta|}$ which  introduces a factor $\frac{1}{|\zeta|}$ in the derivative and the function $\phi$ is replaced by a function having a faster growth close to the boundary of $\Omega$. The operator $\frac{\partial}{\partial x_i}$ is linear continuous from $H_{n,\Phi}$ to  $H_{n+1,\Psi}$ for suitable $\psi$.\\

 Now we will prove that $\mathcal{G}_a(\Omega) \subset \mathcal{G}(\Omega)$. Let us start by  recalling the definition of $\mathcal{G}(\Omega)$:\\
\  \  \\

$ \mathcal{E}(\Omega):=\{f:\Omega
\times ]0,1] \longmapsto \mathbb{C}, \ 
\mathcal{C}^\infty $ such that
\begin{equation} \forall n, \forall K\subset\subset \Omega  
\ \exists N,\ const \ / sup_{x\in K} | \frac{\partial^n f}{\partial x^n}(x,\xi)|  
\leq \frac{const}{\xi^N}   \ 
\forall \xi \ small \ enough \},\end{equation}

$$ \mathcal{N}(\Omega):=\{f \in \mathcal{E}(\Omega) \ / \ \forall n \ \forall K\subset\subset \Omega \ \forall  q \ \exists const>0 \ $$
\begin{equation} / \ sup_{x\in K} |\frac{\partial^n f}{\partial x^n}(x,\xi)| \leq const.\xi^q  \ \forall \xi \ small \ enough \}.\end{equation}
An interesting point is the fact that M. Grosser \cite{GKOS} has proved that   statement (7) is equivalent to the particular case  $n=0$. We define the algebra of nonlinear generalized functions as the quotient\\
$$ \mathcal{G}(\Omega):=\frac{\mathcal{E}(\Omega)}{\mathcal{N}(\Omega)}.$$

It follows obviously from the definition of $H_{n,\phi}$  that if $f\in \mathcal{G}_a(\Omega)$ then the map $(x,\xi)\longmapsto f(x,\xi)$, still denoted by $f$, is in $\mathcal{E}(\Omega)$. This defines a map $i$

$$i: \mathcal{G}_a(\Omega)\longmapsto\mathcal{G}(\Omega),$$ 
$$f(z,\zeta) \longmapsto f(x,\xi)+\mathcal{N}(\Omega),$$
 obtained as the restriction of $f$ to the variables $(x,\xi)$, modulo $\mathcal{N}(\Omega)$.\\

\textbf{Theorem 1.} \textit{The map $i$ is injective, i.e. we can set  $\mathcal{G}_a(\Omega)\subset\mathcal{G}(\Omega)$  to simplify the notation.}\\
\\
\textit{proof}. Assume that $i(f)\in \mathcal{N}(\Omega)$. Therefore
 \begin{equation}\forall K\subset\subset \mathbb{R}, \ \forall  q, \exists  c_q>0 \ /  \ sup_{x\in K} | f(x,\xi)| \leq \ c_q\xi^q \end{equation}
\\
 for $\xi>0$ small enough. Since $f\in H_{n,\phi}$ for some $n,\phi$, it follows from (2,3,4) that for fixed  $x\in O_n$ the map $\zeta \longmapsto f(x,\zeta)$ is holomorphic on the set $0<|\zeta|<\frac{1}{n}$. From the bound (4) this map has a pole at $\zeta=0$, therefore it admits a Laurent series expansion. From (8) all coefficients of the Laurent expansion are null.  Therefore $f(x,\zeta)=0 \ \forall \zeta $ if $x\in O_n$. Fix $\zeta, |\zeta|<\frac{1}{n}$: then  $ f(z,\zeta)=0$ as soon as $x\in O_n, |y_i|<\frac{1}{n} |\zeta|$, because it is null if $z=x\in \mathbb{R}$. Finally since $f$ is 
holomorphic and $V_n$ connected, $f=0$  on $V_n$.$\Box$\\

It is easy to prove (see below)  that this algebra $\mathcal{G}_a(\Omega)$ contains most distributions and it is obvious it contains the algebra of all analytic functions on $\Omega$ whose radius of convergence $r(x)$ has the property $r(x)\geq\alpha>0  \ \forall x$ for some $\alpha>0$ as a faifthful subalgebra : set $F(z,\zeta)=f(z)$ if $f$ is such an analytic function. In order that $\mathcal{G}_a(\Omega)$ would contain all analytic functions on $\Omega$, whose radius of convergence can tend to 0 at infinity or on the boundary of $\Omega$, one should introduce in the definition of $\mathcal{G}_a(\Omega)$ dependence of the bounds in (1,2,4) on the compact subsets of $\Omega$.  The other usual products are recovered through the association \cite{AraBia} p. 105, \cite{Colombeau1} p. 68, \cite{Colombeau2} p. 64, \cite{Ober} p. 93, \dots   only: the algebras of  $\mathcal{C}^\infty$ functions such as the algebra of all $\mathcal{C}^\infty$ functions with compact support are not faithful subalgebras from the use of a convolution to transform them into  holomorphic functions. The assumption $|y_i|\leq\frac{1}{n}|\zeta|$ is used to embed functions and distributions in $\mathcal{G}_a(\Omega)$ by means of a mollifier, see section 5 below.  \\

\textit{Problem}. Is the locally convex inductive limit of the Banach spaces $ H_{n,\phi}$ Hausdorff? We are going to prove a positive answer in the particular case we consider only a suitable countable subfamily of Banach spaces $ H_{n,\phi}$.\\

Now we give two properties of the algebra  $\mathcal{G}_a(\Omega)$  that are quite different from the corresponding properties of  $\mathcal{G}(\Omega)$. Let again $i$ denote the above map:  $\mathcal{G}_a(\Omega)\longmapsto  \mathcal{G}(\Omega)$ and let $f\in \mathcal{G}_a(\Omega)$. Let $\overline{\mathbb{C}}$ denote the ring of generalized complex numbers, i.e. the elements of $\mathcal{G}(\Omega)$ which are constant generalized functions (in the $x$-variable).\\

\textbf{Proposition 1}.  \textit{If all pointvalues $(i(f))(x)\in \overline{\mathbb{C}}$ are null then $f=0$ in  $\mathcal{G}_a(\Omega)$.}\\
\\
\textit{proof}. If $f\in H_{n,\phi}$ and if  $x_0\in O_n$ the pointvalues $(i(f))(x_0)$ have as representatives the analytic functions $\xi\longmapsto f(x_0,\xi)$ which are restrictions of holomorphic functions of $\zeta$ in $\{0<|\zeta|<r\}$ for some $r>0$. The assumption $|f(x_0,\xi)|=o(\xi^q) \  \forall q\in\mathbb{N}$ together with the Laurent series development around the pole $\zeta=0$ (which is a pole from the bound (4)) imply that $f(x_0,\zeta)=0$ in $\{0<|\zeta|<r\}$. Then for fixed $\zeta_0$, consider the holomorphic map $z\longmapsto f(z, \zeta_0)$ which is null if $z=x\in O_n$, as in the proof of theorem 1.$\Box.$\\
\\ 
\textbf{Proposition 2}. \textit{If $\int (i(f))(x) \phi(x) dx=0$ in $\overline{\mathbb{C}}$ for all $\phi\in \mathcal{C}_c^\infty(\Omega)$ then $f=0$ in $\mathcal{G}_a(\Omega)$.}\\
\\
\textit{proof}. If $f\in H_{n,\phi}$ and if $supp(\phi) \subset  O_n$ then \  $\int f(x,\xi) \phi(x) dx=o(\xi^q) \  \forall q$ while the function $\zeta \longmapsto \int f(x,\zeta) \phi(x) dx$ is holomorphic on a set $\{0<|\zeta|<r\}$ with 0 as a pole. As in the proof of theorem 1 the Laurent series expansion permits to conclude that this map is null. Since this holds for all $\phi \in \mathcal{C}_c^\infty(O_n)$ one has $f(x,\zeta)=0$ for all $x\in  O_n$ and all $\zeta \in \{0<|\zeta|<r\}$. One finishes as in Proposition 1.$\Box$\\

%\end{document}

\textbf{3.  Compactness }. \\
\  \  \\
\textbf{Theorem 2}. \textit{$\forall n,\phi \ \ \exists \psi $ such that the canonical inclusion map $$H_{n,\phi} \  \  \  \longmapsto \  \  \   H_{n+1,\psi}$$
$$ \  \  \  \  f:V_n\longmapsto\mathbb{C} \ \  \ \  \  \  \  \  f|_{V_{n+1}}:V_{n+1}\longmapsto\mathbb{C}$$
is a compact operator.}\\  % \end{document}
\\
\textit{proof}.  We seek a function $\psi:\Omega\longmapsto\mathbb{R}^*_+$ depending only on $n$ and $\phi$ such that if $(f_p)_p$ is a sequence in the  unit ball of $H_{n,\phi} $ then we can extract a subsequence that converges in the normed space $H_{n+1,\psi}$ to some limit $f$  in  $H_{n+1,\psi}$. From (4,5)
\begin{equation} \forall p, \ \forall (z,\zeta)\in V_n, \ |f_p(z,\zeta)| \leq \frac{\phi(x)}{|\zeta|^n}.\end{equation}
Therefore the sequence $(f_p)_p$ is uniformly bounded on each compact subset of $V_n$ since $\{\zeta=0\}$ lies on the boundary of $V_n$. We can extract a subsequence,  denoted by $(f_{p(q)})_q$, such that there is a holomorphic function $f$ on $V_n$ such that $f_{p(q)}\rightarrow f$ uniformly on all compact subsets of $V_n$. We intend to find a function $\psi$ such that $ \|f_{p(q)}-f\|_{ H_{n+1,\psi}}\longmapsto 0$, i.e. from (5):
\begin{equation} \forall \epsilon>0 \ \exists q_0 / q\geq q_0 \Rightarrow |(f_{p(q)}-f)(z,\zeta)|\leq\frac{\epsilon \psi(x)}{|\zeta|^{n+1}} \ \forall  (z,\zeta)\in V_{n+1}.\end{equation}

Let $\epsilon>0 $ be given.\\

 $  |(f_{p(q)}-f)(z,\zeta)|\leq |(f_{p(q)})(z,\zeta)|+ |(f)(z,\zeta)|\leq \frac{2\phi(x)}{|\zeta|^n} $ from (9) if $(z,\zeta)\in V_n$.\\
 
 * If $|\zeta|\leq \frac{\epsilon}{2}$, 
\begin{equation}  |(f_{p(q)}-f)(z,\zeta)|\leq 2 |\zeta|\frac{\phi(x)}{|\zeta|^{n+1}} \leq \epsilon \frac{\phi(x)}{|\zeta|^{n+1}}\ \forall  (z,\zeta)\in V_{n+1},|\zeta|\leq \frac{\epsilon}{2} \end{equation}
since $V_{n+1}\subset V_n$.\\

 * If $|\zeta|\geq \frac{\epsilon}{2}, $ \  \ let $K$ be a given compact subset of $\Omega$. The set $K':=\{(z,\zeta)\in V_{n+1} /  \ x\in K, |\zeta|\geq\frac{\epsilon}{2}\}$ is relatively compact in $V_n$: notice one has bounds $\frac{1}{n+1}$ instead of $\frac{1}{n}$ in (1,2). Therefore $f_{p(q)}\rightarrow f$ uniformly on $K'$. Let $|\phi|_K:=sup_{x\in K}|\phi(x)|.$ From this uniform convergence 
\begin{equation} \exists q_0 \ / \ q\geq q_0 \ \Rightarrow |(f_{p(q)}-f)(z,\zeta)| \leq \epsilon |\phi|_K(n+1)^{n+1} \leq  \epsilon\frac{ |\phi|_K}{|\zeta|^{n+1}} \  \forall (z,\zeta)\in V_{n+1}, x\in K, |\zeta|\geq \frac{\epsilon}{2} \end{equation} 
since $|\zeta|<\frac{1}{n+1}$ from (1,2).\\ 

From (11,12) we have obtained that for every compact subset $K$ of $\Omega$ and for every given $\epsilon>0$
\begin{equation} \exists q_0=q_0(K,\epsilon) \ / \ q\geq q_0 \ \Rightarrow |(f_{p(q)}-f)(z,\zeta)| \leq \epsilon \frac{|\phi|_K}{|\zeta|^{n+1}} \ \forall (z,\zeta)\in V_{n+1}, x\in K.\end{equation} 

We have to prove (10), i.e.  the existence of a $q_0$ independent of $K$, with the help that $\phi$ can be replaced by another function $\psi$.\\ 
\\
$\bullet$\textit{Fixed $\epsilon>0$  and variable $K$.}\\
Let $(K_r)_{r\in \mathbb{N}}$ be an exhaustive sequence of compact sets in $\Omega$, i.e. $\Omega=\cup_{r\in \mathbb{N}}K_r$, each compact set in $\Omega$ contained in some $K_r$. The bound (9) implies (since $V_{n+1}\subset V_n$):\\   %\end{document}
$\forall q\in \mathbb{N} \  |(f_{p(q)}-f)(z,\zeta)| \leq 2 \frac{\phi(x)}{|\zeta|^n} \leq  \frac{\phi(x)}{|\zeta|^{n+1}} \forall (z,\zeta)\in V_{n+1}$  (  $|\zeta|<\frac{1}{2}$ if $(z,\zeta)\in V_2)$. Set $\mu_r:=\frac{|\phi|_{K_r}}{\epsilon}.$ Then if $ (z,\zeta)\in V_{n+1}$ with $x\in K_r$ we have \\
$\forall q\in \mathbb{N} \ |(f_{p(q)}-f)(z,\zeta)| \leq \epsilon \frac{\mu_r}{|\zeta|^{n+1}}.$ 
 We have proved that for all given $\epsilon>0$ there exists a sequence $(\mu_r)_r$ such that 
\begin{equation} \forall q\in \mathbb{N}  \ \ |(f_{p(q)}-f)(z,\zeta)| \leq \epsilon \frac{\mu_r}{|\zeta|^{n+1}} \ \forall (z,\zeta)\in V_{n+1}, x\in K_r.\end{equation}
 The sequence $(\mu_r)_r$ we have constructed depends on $\epsilon>0$. We set $\mu_{K_r,\epsilon}$ for $\mu_r$ in (14). Now we will find a sequence  $(\nu_r)_r$,   independent of $\epsilon>0$, using (13) and a standard manipulation of sequences.\\ 
 \\
$\bullet$ \textit {Let $\epsilon\rightarrow 0$.}
For this let us consider a sequence $(\epsilon_n)_n$, $\epsilon_{n+1}<\epsilon_n \forall n, \epsilon_n\rightarrow 0$.\\

* If $\epsilon=\epsilon_1 $ (14) gives a sequence  $(\mu_r)_r$, that we denote $(\mu_{K_r,\epsilon_1})_r$, such that \\
$q\geq q_0(\epsilon_1) \Rightarrow  |(f_{p(q)}-f)(z,\zeta)| \leq \epsilon_1 \frac{\mu_{K_r,\epsilon_1}}{|\zeta|^{n+1}} \ \forall  (z,\zeta)\in V_{n+1}, x\in K_r$,(with $q_0(\epsilon_1)=1$ because indeed from (14) this holds for all $q$).\\
Choose a sequence  $(\nu_r)_r$, to be modified later inductively  by letting fixed the first $n$ terms at the step n of the induction, such that 
\begin{equation} \nu_r\geq \mu_{K_r,\epsilon_1} \ \forall r\geq 1.\end{equation}
Then  
\begin{equation} q\geq q_0(\epsilon_1)\Rightarrow  |(f_{p(q)}-f)(z,\zeta)| \leq \epsilon_1 \frac{\nu_r}{|\zeta|^{n+1}} \ \forall  (z,\zeta)\in V_{n+1}, x\in K_r.\end{equation}

* If $\epsilon=\epsilon_2 $  \  keep $\nu_1$ unchanged  and increase the values $\nu_2,\nu_3, \dots$ such that $\nu_r\geq\mu_{K_r,\epsilon_2}$ if $r\geq 2$. Then we claim that $\exists q_0(\epsilon_2)$ such that
\begin{equation} q\geq q_0(\epsilon_2)\Rightarrow  |(f_{p(q)}-f)(z,\zeta)| \leq \epsilon_2 \frac{\nu_r}{|\zeta|^{n+1}} \ \forall  (z,\zeta)\in V_{n+1}, x\in K_r.\end{equation}
Indeed: when $r\geq 2$ this follows from (14) and holds for all $q$; when $r=1$
use (13) with $K=K_1$, and  $\epsilon$ such that $\epsilon |\Phi|_{K_1} \leq \epsilon_2 \nu_1$; then $q_0(\epsilon_2)=q_0(K_1,\epsilon) $ from (13).\\

* If $\epsilon=\epsilon_3 $ \ keep $\nu_1$ and $\nu_2$ unchanged  and increase the values $\nu_3,\nu_4, \dots$ such that $\nu_r\geq\mu_{K_r,\epsilon_3}$ if $r\geq 3$. Then we claim that
\begin{equation} q\geq q_0(\epsilon_3)\Rightarrow  |(f_{p(q)}-f)(z,\zeta)| \leq \epsilon_3 \frac{\nu_r}{|\zeta|^{n+1}} \ \forall  (z,\zeta)\in V_{n+1}, x\in K_r.\end{equation}
Indeed: when $r\geq 3$ this follows from (14) and it holds for all $q$; when $r=1,2$ use (13) with $K=K_1,K_2$ and $\epsilon$ such that $\epsilon |\Phi|_{K_1} \leq \epsilon_3 \nu_1$ and $\epsilon |\Phi|_{K_2} \leq \epsilon_3 \nu_2$;  then $q_0(\epsilon_3) \geq max(q_0(K_1,\epsilon),q_0(K_2,\epsilon)) $ from (13). \\

 The induction is obvious. Finally we have obtained:
\begin{equation} \forall \epsilon>0 \  \exists q_0(\epsilon) \ / \ q\geq  q_0(\epsilon)\Rightarrow  |(f_{p(q)}-f)(z,\zeta)| \leq \epsilon \frac{\nu_r}{|\zeta|^{n+1}} \ \forall  (z,\zeta)\in V_{n+1}, x\in K_r.\end{equation}

Now it suffices to choose a function $\psi$ such that $\forall r \  \psi(x)\geq \nu_r$ if $x\in K_r, x\not\in K_{r-1}$ to have that 
\begin{equation} \forall \epsilon>0  \ \exists q_0(\epsilon) \ / \ q\geq  q_0(\epsilon)\Rightarrow  |(f_{p(q)}-f)(z,\zeta)| \leq \epsilon \frac{\psi(x)}{|\zeta|^{n+1}} \ \forall  (z,\zeta)\in V_{n+1} \end{equation}
$\Box$

%\end{document}

\textbf{4. Construction of a Hausdorff locally convex subalgebra of $\mathcal{G}(\Omega)$.}
In the algebra  $\mathcal{G}_a(\Omega)$ we select a countable increasing sequence  of Banach  spaces $(H_{n_p,\phi_p})_p$  that satisfy the following three properties:
\\
\\ 
i) $H_{n_p,\phi_p}. H_{n_p,\phi_p} \subset H_{n_{p+1},\phi_{p+1}},$ with continuous inclusion,\\ 
ii) $\frac{d}{dx}H_{n_p,\phi_p}\subset H_{n_{p+1},\phi_{p+1}},$ with continuous inclusion,\\
iii) the inclusion map $ H_{n_p,\phi_p} \longmapsto  H_{n_{p+1},\phi_{p+1}}$ is a compact operator.\\
         
 Then it is known \cite{Bourbaki} Prop.7, EVT III.6 ( see also \cite{Grothendieck}, ex. 4 p. 145) that the locally convex topological inductive limit of the normed spaces $ H_{n_p,\phi_p}$ is Hausdorff, and that any bounded set in this locally convex space is contained and bounded in some normed space $ H_{n_p,\phi_p}$. Further it is  known  from the theory of locally convex topological vector spaces that such an inductive limit  space is a strong dual of a  Fr\'echet-Schwartz space (DFS space for short):  a proof can be found in the books \cite{Perez} p. 293, \cite{Hogbe2} p. 99.\\
 
To prove that the locally convex inductive limit of the normed spaces $H_{n_p,\phi_p}$ is a topological algebra, i.e. that the multiplication is continuous, one has to prove that for every bornivorous disk $Q$ there is a bornivorous disk $Q'$ such that $Q'.Q'\subset Q$ since the bornivorous disks form a fundamental system of 0-neighborhoods. For this proof let $Q$ be a given bornivorous disk and denote by $B_n$  balls in each of the normed spaces $H_{n_p,\phi_p}$, that we can choose such that $B_n.B_n\subset B_{n+1}$ for simplicity. If $S$ is a set in a vector space we denote by $\Gamma(S)=\{\sum_{finite}\lambda_i x_i\}_{\sum|\lambda_i|\leq 1, x_i\in S}$ 
the convex  balanced  span of $S$. It is easy to prove by induction the existence of a sequence $(\epsilon_n>0)_n$ such that $\Gamma((\sum_{i=1}^n \epsilon_iB_i).(\sum_{i=1}^n \epsilon_iB_i))\subset(1-\frac{1}{n+1})Q$. If $(A_i) $ are subsets of a vector space set $\sum_{i=1}^\infty A_i:=\{\sum_{finite}
 x_i\}_{x_i\in A_i}$. Then we have obtained that  $\Gamma((\sum_{i=1}^\infty\epsilon_iB_i).(\sum_{i=1}^\infty \epsilon_iB_i))\subset Q$. Now if $Q':= \sum_{i=1}^\infty\epsilon_iB_i$ one has $Q'.Q'\subset Q$.\\

 %\end{document}

We denote any such locally convex topological inductive limit by $\mathcal{A}(\Omega)$. It is a  subalgebra of $\mathcal{G}(\Omega)$ having the topological and compactness   properties of DFS spaces, in particular any bounded set is contained and bounded in the normed space spanned by a bounded disk, normed with the gauge of this bounded disk. The DFS algebra $\mathcal{A}(\Omega)$ put in evidence from this construction depends on the choice of the sets $(O_n)_n$ and then on the choice of the increasing sequence of Banach spaces $(H_{n_p,\phi_p})_p$ with compact inclusions.\\

\textbf{5. Inclusion    of distributions.}
 We give the proof in the case $\Omega=\mathbb{R}$ and with sets $O_n=]-\infty,n[\cup]n,+\infty[$ or $]-\infty,n[$, or $]n,+\infty[$. The proof in the general case is similar. We use the standard analytic mollifier  $ \rho(z)=const.\frac{1}{(1+z^2)^s}, s\in \mathbb{N}$ large enough so as to permit integration, as needed in $\mathbb{R}^n, n>1$ (easily observed from the use of spherical coordinates).\\
\\
 %\end{document}
\textbf{Proposition 3}: \textit{The map  $ j:f\longmapsto F$ defined by $F(z,\zeta)=(f*\rho_\zeta)(z),\    \rho_\zeta(z)=\frac{1}{\zeta}\rho(\frac{z}{ \zeta}),\ \rho(z)=\frac{1}{\pi}.\frac{1}{1+z^2}, \ \int\rho(x)dx=1,$  maps the vector space $\mathcal{E}'(\mathbb{R})$ of all distributions with compact support into an algebra  $\mathcal{A}(\mathbb{R})$. The map $j$ is injective and continuous from $\mathcal{E}'(\mathbb{R})$ into $\mathcal{A}(\mathbb{R})$.}\\

% \end{document}
proof. From the classical structure theorem of distributions with compact support it suffices to consider the embedding of a continuous function with compact support. If $ supp(f)\subset [-A,+A]$ then
 $$F(z,\zeta)=\frac{1}{\pi}.\zeta^{-1}\int_{-A}^{+A}f(\lambda)\frac{1}{1+(\frac{\lambda-z}{\zeta})^2}d\lambda.$$
%\end{document}
The function to be integrated on $[-A,A]$ is defined and continuous as long as $\zeta\not=0$ (which holds in all sets $V_n$) and 
\begin{equation}   \frac{\lambda-z}{\zeta}\not=\pm i  \  \ \forall \lambda \in [-A,A]. \end{equation}
%\end{document}
Therefore in the open set in $(z,\zeta)$ for which these conditions are satisfied the function $F$ is well defined and holomorphic. (21) reads as 
\begin{equation}    \lambda-x\pm \eta \not=\pm i(\xi\pm y). \end{equation}
%\end{document}
We distinguish two cases.\\ 
\\
$\bullet$If $|x|>A+r$ for some $r>0$ then $|\lambda-x|>r \ \forall \lambda \in [-A,A].$ Therefore if $0<|\zeta|<r$, hence $|\eta|<r$, the left member of (22) is real nonzero, therefore (22) holds. \\
\\
$\bullet$Now if $|x|\leq A+r$ then (22) is ensured if $\xi>0$ and $|y|<\xi$ since  then the right member is nonzero purely imaginary.\\

 Setting here $O_n=]-\infty,-n[\cup]n,+\infty[$ one obtains at once from (1,2) that $F$ is defined and holomorphic on $V_n$ for $n$ large enough. On such a set $V_n$ 
 there is $\alpha>0$ such that  $|\frac{\lambda-z}{\zeta}\pm i|\geq \alpha >0 \  \forall \lambda \in [-A,A]$. Therefore  the function $\frac{1}{1+(\frac{\lambda-z}{\zeta})^2}$ is uniformly bounded if  $\lambda$ ranges in $[-A,A]$. 
Therefore 
%\end{document}
\begin{equation} |F(z,\zeta)| \leq \frac{\frac{1}{\pi}}{|\zeta|} \  \forall (z,\zeta) \in V_n \end{equation}
$\Box$.

\  \ \\   
\textit{ 
This result can be extended to distributions constant at infinity as follows.} Let $f$ be an $L^p_{loc}$ function constant at infinity. Then if $c$ denotes the value of $f$ near $\infty$ \\ 
$$F(z,\lambda)=\frac{1}{\pi}.\zeta^{-1}(\int_{-\infty}^{+\infty}(f(\lambda)-c) \frac{1}{1+(\frac{\lambda-z}{\zeta})^2} d\lambda+\int_{-\infty}^{+\infty}c.\frac{1}{1+(\frac{\lambda-z}{\zeta})^2}   d\lambda ).$$
 $f(\lambda)-c$ has compact support therefore the proof of  proposition 3 applies to the first integral. Concerning the second integral one has (from the choice of $const$ so as to have a mollifier):
$$\frac{1}{\pi}.\zeta^{-1}\int_{-\infty}^{+\infty}\frac{1}{1+(\frac{\lambda-z}{\zeta})^2} 
d\lambda=1.$$
Indeed, first, for fixed $x$ the holomorphic function $\zeta \longmapsto \frac{1}{\pi} .\zeta^{-1}  \int_{-\infty}^{+\infty} \frac{1}{1+(\frac{\lambda-x}{\zeta})^2}  d\lambda$ is  constant $=1$ since, by definition of  a mollifier, it is $=1$ for  $\zeta =\xi$ and it is holomorphic in $\zeta$; then  fix $\zeta$ and let $z$ vary: one obtains  from holomorphy that  $z \longmapsto \frac{1}{\pi} .\zeta^{-1}\int_{-\infty}^{+\infty} \frac{1}{1+(\frac{\lambda-z}{\zeta})^2}  d\lambda=1$ since it is for $z=x$. Therefore 
$$F(z,\zeta)=c+\frac{1}{\pi}.\zeta^{-1}\int_{-\infty}^{+\infty}(f(\lambda)-c) \frac{1}{1+(\frac{\lambda-z}{\zeta})^2} d\lambda
$$  is defined on some set $V_n$ for $n$ large enough, holomorphic there and bounded of the form (23) there. This extends obviously to distributions   by derivation. If the presentation is done with $O_n$ sets only located at $+\infty$ i.e. of the form $]-\infty,-A[ $ or  $]A,+\infty[, A>0$,  then $f$ can have different constant values at $+\infty$ and $-\infty$. If the presentation is done with $O_n$-type sets located around a real number then these results can be extended to distributions with a rather fast growth at infinity by choosing mollifiers $const. \frac{1}{(1+(\frac{\lambda-z}{\zeta})^2)^s}$ with $s$ large enough to ensure mollifiers having a fast enough decrease at infinity. By derivation and integration this result can be extended to distributions that are polynomials at $\infty$. \\

%\end{document}

\textbf{ 6. Connection with the Scarpalezos' sharp topology.} In this section we prove a coherence result between the topological convergence  considered above and convergence in the sense of the sharp topology: when restricted to bounded subsets of $\mathcal{G}_a(\Omega)$  the sharp topology (which basically deals with "infinitesimals") is strictly finer than the classical  topology defined in this paper.\\

A fundamental system of 0-neighborhoods for the sharp topology in $\mathcal{G}(\Omega)$  is made of the sets \\

$ V(K,p,q)=\{G\in \mathcal{G}(\Omega) $ such that there exists a  representative    R  of  G  such that
\begin{equation} \exists \eta>0,C>0/ |DR(x,\xi)|\leq C\xi^q \ if \ order(D) \ \leq p ,x\in K, 0<\xi<\eta.\} \end{equation}
Then the same statement holds for any representative (with different $\eta,C$). The neighborhoods of any point are obtained  by translation. We recall the sharp topology is Hausdorff, metrizable and complete but that it is not a vector space topology \cite{AragonaFernandez, AragonaFerJu3, Scarpatopo, Garetto, Garetto2}.\\

$\mathcal{G}_a(\Omega)$ can be considered as an algebraic inductive limit of the Banach spaces $H_{n,\phi}$ and we naturally say that a subset of $\mathcal{G}_a(\Omega)$  is bounded iff it is contained and bounded in one of these Banach spaces. A sequence is said to be convergent iff it is convergent in one of these Banach spaces. This particular convergence structure has been studied by various authors (Waelbroeck, Sebasti\~ao e Silva, Hogbe-Nlend, \dots), see \cite{Colombeau0} for more details. This convergence can be related to the convergence for the sharp topology as follows. Let $i$ denote the inclusion of $\mathcal{G}_a(\Omega)$ into $\mathcal{G}(\Omega)$ considered in theorem 1.\\

\textbf{Proposition 4. }\textit{ If $(f_{n})_{ n\in \mathbb{N}}$ is  bounded in  $\mathcal{G}_a(\Omega)$ and if $f\in \mathcal{G}_a(\Omega)$, then $i(f_n)\rightarrow i(f)$ for the sharp topology in  $\mathcal{G}(\Omega)$ implies that $f_n\rightarrow  f$ in  $\mathcal{G}_a(\Omega)$.}\\

In other words the sharp topology induces on the bounded sets of  $\mathcal{G}_a(\Omega)$ a topology finer than the one induced by the convergence structure of  $\mathcal{G}_a(\Omega)$. The converse is wrong: if $f\in  \mathcal{G}_a(\Omega)$ the sequence $(\frac{1}{n} f)_n$ tends to $0$ in  $\mathcal{G}_a(\Omega)$ and  $i(\frac{1}{n} f)$ does     not tend to $0$ for the sharp topology in  $\mathcal{G}_a(\Omega)$ if $f\not= 0$.\\

\textit{proof of Prop.4.}  From theorem 2 (compactness) we can assume, modulo extraction of a subsequence, that $(f_n)_n$ converges (in some suitable Banach space of the algebraic inductive limit) to some $g\in \mathcal{G}_a(\Omega)$.  The letter $K$ denotes a compact set. It suffices  to prove that $g=f$: then any convergent subsequence of $(f_n)_n$ would also converge to $f$ and proposition 4 would be proved. The assumption implies that \\
$$ \forall K,\forall q \  \exists n_0 / n\geq n_0 \Rightarrow i(f_n)-i(f)\in V(K,0,q).$$
This can be rewritten as
\begin{equation} \forall K,q \  \exists n_0 / n\geq n_0 \Rightarrow \exists C_n,\eta_n>0 / |(f_n-f)(x,\xi)|<C_n\xi^q \  if \ 0<\xi<\eta_n, x\in K.\end{equation}
The difficulty in the proof comes from the fact there is no control on $C_n$ and $\eta_n$, which therefore cannot be used inside a bound.
Since the set $\{f_n, n\in \mathbb{N}\}\cup \{f\} $ is a bounded set in $ \mathcal{G}_a(\Omega)$ it is contained and bounded in some Banach space $H_{p_0,\phi_0}$. Therefore if one chooses $K\subset O_{p_0}$ there exists a fixed value $R$ independent of $n$ such that all $f_n-f$ are holomorphic functions of $\zeta, 0<|\zeta|<R$, and are
uniformly bounded there by $\frac{M}{|\zeta|^{p_0}}$, where $M$ is some fixed constant, for any fixed $x\in K$. \\
%\end{document}
The point $\zeta=0$ which a priori is a pole of $f_n-f$ (from the bounds in definition of $ \mathcal{G}_a(\Omega)$) is a removable singularity from (25), for $n$ large enough. Therefore each function $\zeta \longmapsto (f_n-f)(x,\zeta)$ can be developped in a  Taylor series at $\zeta=0$ in the disk $|\zeta|<R$. Cauchy 's inequalities give 
$$ \zeta^{p_0}(f_n-f)(x,\zeta)=\sum_{i=0}^\infty a_i \zeta^i$$ where the coefficients $a_i$ depend on $n$ but satisfy the uniform bounds $|a_i|\leq \frac{M}{R^i}$ for all $n$. From (25) the coefficients $a_i$ of $f_n-f$ are null if $i<q-p_0$ for $n\geq n_0$. Therefore the uniform bound implies that for $n\geq n_0$ \\
$$|\zeta^{p_0}(f_n-f)(x,\zeta)|\leq M.|\frac{\zeta}{R}|^q (1+|\frac{\zeta}{R}|+|\frac{\zeta}{R}|^2+\dots)\leq const.|\frac{\zeta}{R}|^q$$
if $|\zeta|\leq \frac{R}{2}$.\\

It suffices now to use the convergence of $f_n$ to $g$ in the Banach space  $H_{p_0,\phi_0}$: \\
$$|(f-g)(x,\zeta)| \leq |(f_n-f)(x,\zeta)|+|(f_n-g)(x
,\zeta)| \leq const.|\frac{\zeta}{R}|^q\frac{1}{|\zeta|^{p_0}}+\|f_n-g\|_{H_{p_0,\phi_0}}\frac{\|\phi_0\|_K}{|\zeta|^{p_0}}$$
and let $q$ tend to $\infty$ with $n\geq n_0=n_0(q)$ (from (25); we used (25) only to state that $f_n-f$ is of order $\geq q$ at 0 since we have no control on $C_n,\eta_n$). Therefore $g=f$ on K. One can choose $K$ having a nonvoid interior and finish as in the proof of theorem 1 using uniqueness of analytic continuation. $\Box$\\

%\end{document}

\textbf{7. Nuclearity.}  In this section we prove that, with $O_n=\{x\in \mathbb{R}^{n_0} / |x|>n\}$, the space $\mathcal{E}'$ of all distributions on $\mathbb{R}^{n_0}$ with compact support is contained in a  differential algebra $\mathcal{A}(\mathbb{R}^{n_0})$ which is a strong dual of a nuclear Fr\'echet space (DFN space for short). By definition  the space $\mathcal{E}'$ is the strong dual of the nuclear Fr\'echet space   of all  $\mathcal{C}^\infty$  functions on $ \mathbb{R}^{n_0}$ with its classical topology of uniform convergence on compact sets of the functions and all partial derivatives. This will show that the DFS spaces constructed in section 4 could be replaced by the considerably richer structure of DFN spaces.\\

%\end{document}
\textit{Recall of a classical notation.} If a sequence $(x_n)$ tends to 0 in a Banach space $E$ we set%\end{document}
$$ \Gamma_{l_1} (x_n)=\{ \sum_n \lambda_n  x_n\}_{\sum |\lambda_n| \leq 1}$$ 
\\
which is  the closed convex balanced hull of the sequence $(x_n)$ and is the same for any larger Banach space $E$ (one proves it is compact in $E$). It is obvious that, in a Banach algebra $E$ or more generally in an "algebra stemming from an algebraic inductive limit of Banach spaces" like $\mathcal{G}_a$,  the product of a finite number of $ \Gamma_{l_1} $ sets is still such a  set: $ \Gamma_{l_1} (x_n). \Gamma_{l_1} (y_p)=\{ \sum_n \lambda_n  x_n\}_{\sum |\lambda_n| \leq 1}.\{ \sum_p \mu_p  y_p\}_{\sum |\mu_p| \leq 1}=\{ \sum_{n,p} \lambda_n \mu_p  x_n.y_p\}_{\sum_{n,p} |\lambda_n.\mu_p| \leq 1}$ and it suffices to order the double sequence $x_n.y_p$ into a simple sequence.\\

\textit{Recall of a classical result.} In the $\mathcal{G}_a$-context it is particularly convenient to use the classical results on nuclear spaces in the way they are formulated in \cite{Hogbe1}. In particular we will use a theorem  of \cite{Hogbe1} p. 82   which is a characterisation of algebraic inductive limits of Banach spaces (called bornological vector spaces there but this terminology can be misleading) that have the nuclearity property in this context. An algebraic inductive limit $E$ of Banach spaces (such as  $\mathcal{G}_a(\Omega)$) is said to be nuclear iff for any of these Banach spaces $E_i$ there is another of these Banach spaces, say $E_j$, which contains $E_i$ such that the inclusion from $E_i$ into $E_j$ is a nuclear map. Then it follows from classical results on nuclear maps, \cite{Hogbe1} p. 82,  that  $E$ is such a nuclear algebraic  inductive limit of Banach spaces iff any bounded set (i.e. by definition a set which is contained and bounded  in a Banach space from the inductive limit) is contained in a  $\Gamma_{l_1} (x_n)$ set for a sequence $(x_n)$ that tends to 0 in the normed space spanned by a bounded set. After these recalls we will prove:\\

\textbf{Proposition 5.} \textit{The locally convex vector space $\mathcal{E}'$ of all distributions on $\mathbb{R}^n$ with compact support is continuously contained (through the map $j$ of proposition 3) in a differential algebra $\mathcal{A}(\mathbb{R}^n)$ which is a strong dual of a nuclear Fr\'echet space}.\\ 

In short the classical DFN space $\mathcal{E}'$, which is not an algebra, is embedded into a  DFN algebra.\\

\textit{Sketch of proof}. The technical proof is somewhat sketched to make it more accessible. The DFN space $\mathcal{E}'$  is equipped with its usual bounded sets: as a DFN space it has a fundamental sequence $(b_n)$ of bounded disks such that for all n there exists a sequence $(x_n^p)_p$ of distributions ($\in\mathcal{E}'$) which tends to 0 in the Banach space $(\mathcal{E}')_{b_{n+1}}$ such that 
$$b_n\subset  \Gamma_{l_1} (x^n_p)_p.$$
This  structure result of DFN spaces follows from  the classical result recalled above. This is transfered to $j(\mathcal{E}')$ if $j:\mathcal{E}'\rightarrow \mathcal{G}_a$ is the obvious extension to distributions of the map $j$ in section 5. In $ \mathcal{G}_a$ we set, with $f^n_p=j(x^n_p)$,
$$j(b_n)\subset  \Gamma_{l_1} (f^n_p)_p.$$%\end{document}
Since the product in $ \mathcal{G}_a$ of two $ \Gamma_{l_1}$  sets is still a  $\Gamma_{l_1}$  set,
$$j(b_1).j(b_1)\subset  \Gamma_{l_1} (F^1_p)$$%\end{document}
 for a sequence $(F^1_p)_p$ in $ \mathcal{G}_a$ which tends to 0 in the Banach space spanned by a bounded disk. There is a sequence of real numbers $\nu^1_p$, which tends to $+\infty$, such that the sequence $(G^1_p:=\nu^1_p F^1_p)_p$ tends again to 0 in the same Banach space as the sequence $(F^1_p)_p$. Set
 $$B_1:= \Gamma_{l_1} (G^1_p)$$%\end{document}
which is a bounded disk in  $ \mathcal{G}_a$ such that $ j(b_1).j(b_1)\subset B_1$. Now we go on the inductive construction of a suitable sequence $(B_n)$ of bounded sets in $ \mathcal{G}_a$. As above since $j(b_2)$ is a   $\Gamma_{l_1}$ set one has %\end{document}
$$(j(b_2))^3.(B_1)^2\subset \Gamma_{l_1}(F_p^2)$$
where $S^2=S.S=\{x.y\}_{x,y\in S}$, etc, and we set as above 
$$B_2 := \Gamma_{l_1} (G^2_p).$$
for some sequence $(G^2_p)_p$. 
The induction is obvious: by considering $(j(b_r))^{r+1}.(B_{r-1})^{r}$ as starting point one constructs a bounded disk
$$B_r := \Gamma_{l_1} (G^r_p)$$%\end{document}
of  $ \mathcal{G}_a$. From the theorem \cite{Hogbe1} recalled above (converse part of the characterization) the inductive limit of the Banach spaces  $ (\mathcal{G}_a)_{B_n}$  is a nuclear inductive limit of Banach spaces as well as it is a subalgebra of  $ \mathcal{G}_a$. Similarly as  the compactness case in section 4 this nuclear inductive limit is a DFN space (instead of a DFS space in section 4) because of the countability of the family $(B_n)$ and the fact that a nuclear map is a fortiori compact. However the proof is not finished: one should also include into this inductive construction all sets made of partial derivatives of the functions in the various  $\Gamma_{l_1}$-sets we have constructed. Since the set of  partial derivatives (up to a given order) of elements of a $ \Gamma_{l_1}$-set is obviously the  $\Gamma_{l_1}$-set of the corresponding partial derivatives, and since all these sets of derivatives that appear in the inductive construction can obviously be inserted into the inductive construction above, finally one can obtain an inductive limit of Banach spaces which is as above and further is closed under all partial derivatives, which proves the Proposition.$\Box$\\

Finally, using the $\Gamma_{l_1}$-sets and the characterization of nuclearity through these sets, the  construction of DFN subalgebras of   $ \mathcal{G}_a$ is similar to the construction of the DFS subalgebras in section 4.\\

%\end{document}

\textbf{8. Extension to a differential calculus.} Let $\Omega$ be a connected open set in $\mathbb{R}^n$ and let $\tilde{\Omega}_c$ be the set of all "`generalized points in $\Omega$"', i.e. the classes in $\overline{\mathbb{R}^n}$ of  maps $[0,1[\longmapsto \mathbb{R}^n$, $\epsilon\longmapsto x(\epsilon)$, which map $[0,1[ $ into a compact set in $\Omega$ and are moderate (concerning the derivatives). $\tilde{\Omega}_c$ is an open set in  $\overline{\mathbb{R}^n}$ for the sharp topology \cite{Kunzinger,AragonaFerJu1}: a systematic study of versions with smooth or continuous dependence on $\epsilon$ has been done in \cite{Burtscher}. The authors of   \cite{AragonaFerJu1}
introduced a faithful  extension of $\mathcal{G}(\Omega)$ in form of an original differential calculus for maps $\tilde{\Omega}_c\longmapsto\overline{K}, K=\mathbb{R}$ or $\mathbb{C}$: using the sharp topology they define $\mathcal{C}^\infty$ maps from $\tilde{\Omega}_c$ into $\overline{K}$, and $\mathcal{G}(\Omega)\subset\mathcal{C}^\infty(\tilde{\Omega}_c,\overline{K})$ as a faithful differential algebra. This last differential algebra appears therefore as a larger differential algebra of nonlinear generalized functions whose presentation is far closer to the usual presentation of differential calculus. In this sction we sketch a similar differential calculus in the case of the differential algebra  $\mathcal{G}_a(\Omega)$. It has been proved in Proposition 1 that  an element $f$ in $\mathcal{G}_a(\Omega)$ is entirely defined by the set of all pointvalues on points in $\Omega$, while an element of $\mathcal{G}(\Omega)$ is only defined by the knowledge of all its pointvalues on the "`generalized points"' elements of $\tilde{\Omega}_c$ \cite{Kunzinger}. Therefore the domain of the analogous differential calculus will now be $\Omega$ itself, with values in some ring of generalized numbers adapted to the  $\mathcal{G}_a$ setting.\\

We define $H_n:=\{\zeta\longmapsto f(\zeta)$ which are holomorphic in the sector $|arg\zeta|<\frac{1}{n}, \ 0<|\zeta|<\frac{1}{n}$ and such that $$ |f(\zeta)|<\frac{const}{|\zeta|^n}$$ there.\}. The vector space $H_n$ is equipped with the norm\\
\begin{equation} \|f\|_n=sup_{|arg\zeta<\frac{1}{n},0<|\zeta|<\frac{1}{n}}|\zeta|^n |f(\zeta)|. \end{equation}
We define the set $\overline{K}_a$ of generalized numbers as the algebraic inductive limit of the normed spaces $A_n$ analogously to the definition of $\mathcal{G}_a$ as an inductive limit.\\

\textbf{Proposition 6.} \textit{The inclusion map $H_n\longmapsto H_{n+1}$ is a compact operator.}\\
\\
The proof is similar to a part of the proof of theorem 1: the $x$-variable is absent. \\ %\end{document}

Therefore   $\overline{K}_a$ is a strong dual of a Fr\'echet-Schwartz space. The ring $\overline{K}_a$ is not a field: the element  $\zeta\longmapsto exp(-\frac{1}{\zeta^2})$ is nonzero and noninvertible. The differential algebra of all $\mathcal{C}^\infty$ maps from $\Omega$ to  $\overline{K}_a$, \cite{Colombeau0} for instance, could play the role of the differential algebra $\mathcal{C}^\infty(\tilde{\Omega}_c,\overline{K})$ defined in \cite{AragonaFerJu1}.\\
\\

%\end{document}

\textbf{9. Conclusion.} In contrast with the original introduction of nonlinear generalized functions in \cite{Colombeauoriginal} and the exposition of this theory in various expository texts \cite{Colombeau1,Colombeau2, Bialivre, AraBia, Christiakov, Ober, Colombeau3, Colombeau4,Hoskins, Nedel, GKOSmemoirs, GKOS, SteinVi, Todorov}, we have put in evidence a very rich hidden locally convex topological structure in subalgebras of $\mathcal{G}(\Omega)$ that permits the use of the classical deep tools of topology and functional analysis.  Many variants of the original  construction introduced in this paper are possible.\\
%\end{document}

One could insist on the difference between this context in which one has both compactness  (from any bounded sequence  in one of these locally convex subalgebras, one  can extract a convergent subsequence), and compatibility with nonlinearity (in these topological subalgebras  $F_n \rightarrow F, G_n \rightarrow G \Rightarrow F_n.G_n \rightarrow F.G$), as well as compatibility with partial derivatives. In contrast, in classical mathematics, there is a well known conflict between nonlinearity and compactness: the strong topology in Banach spaces is often compatible with nonlinearity but not with compactness (Riesz theorem) and the weak topologies provide compactness but are incompatible with nonlinearity.\\

 In some sense this paper starts the extension of the rich topological aspects of Schwartz distribution theory to the nonlinear context:  one can even  benefit from  very deep properties such as those from nuclearity: the whole  theory of nuclear  locally convex spaces and nuclear maps \cite{Pietsch} can now be applied in a nonlinear setting of differential algebras of generalized functions suited for mathematical analysis.\\

Such a topological context  was lacking in the previous expositions of nonlinear generalized functions  because  the explicit presence of a quotient was a nuisance for the topological structure. On the other hand this quotient permits to identify objects that give same results in multiplication of distributions and permits to present the nonlinear theory as a direct (algebraic) extension of distribution theory \cite {Colombeauoriginal, Colombeau1}. In the present paper there is in fact compatibility of this quotient and Hausdorff locally convex topological structures in specific subalgebras of the whole algebra of nonlinear generalized functions in which the quotient disappears in practice from the fact one can single out a unique privilieged representative for the generalized functions in these subalgebras.\\
\\
Acknowledgements. The authors are very indebted to the referee for numerous improvements in the text.

%\end{document}

\end{document}